\documentclass{amsart}
\usepackage{amssymb}
\usepackage{graphicx}
\usepackage[all]{xy}
\usepackage{hyperref}
\usepackage{a4wide}
\newtheorem{lemma}{Lemma}[section]
\newtheorem{theorem}[lemma]{Theorem}

\newtheorem{corollary}[lemma]{Corollary}

\newtheorem{proposition}[lemma]{Proposition}
\theoremstyle{remark}

\newtheorem{remark}[lemma]{Remark}

\newtheorem{problem}[lemma]{Problem}

\theoremstyle{definition}

\newtheorem{definition}[lemma]{Definition}
{}
\renewcommand{\phi}{\varphi}

\def\id{{\rm id}}
\def\ord{{\rm ord}}

\def\normal{\triangleleft}
\newcommand{\gal}{\textnormal{Gal}}
\newcommand{\bfsigma}{\mbox{\boldmath$\sigma$}}
\def\ab{{\rm ab}}
\def\Fhat{{\hat F}}

\def\Ggag{{\bar{G}}}
\def\alphagag{{\bar{\alpha}}}
\def\betagag{{\bar{\beta}}}
\def\thetagag{{\bar{\theta}}}
\def\etagag{{\bar{\eta}}}
\def\Agag{{\bar{A}}}
\def\Bgag{{\bar{B}}}
\def\Hgag{{\bar{H}}}

\def\calS{\mathcal{S}}

\newcommand{\bbC}{\mathbb{C}}
\newcommand{\bbF}{\mathbb{F}}
\newcommand{\bbN}{\mathbb{N}}
\newcommand{\bbQ}{\mathbb{Q}}
\newcommand{\bbZ}{\mathbb{Z}}

\def\alphahat{{\hat \alpha}}
\newcommand{\ihat}{\hat{i}}
\newcommand{\jhat}{\hat{j}}
\newcommand{\muhat}{\hat{\mu}}
\newcommand{\nuhat}{\hat{\nu}}
\newcommand{\betahat}{\hat{\beta}}
\newcommand{\etahat}{\hat{\eta}}
\newcommand{\thetahat}{\hat{\theta}}

\def\Ahat{{\hat A}}
\def\Bhat{{\hat B}}
\def\Ghat{{\hat G}}
\def\Hhat{{\hat H}}
\def\Nhat{{\hat N}}
\def\Mhat{{\hat M}}
\def\thetahat{{\hat \theta}}

\newcommand{\bfb}{\mbox{{\rm \textbf{b}}}}
\newcommand{\bfc}{\mbox{{\rm \textbf{c}}}}

\newcommand{\bfh}{\mbox{{\rm \textbf{h}}}}
\newcommand{\invlim}{{\displaystyle \lim_{\longleftarrow}}}

\newcommand{\C}{\mathcal{C}}

\begin{document}

\fontsize{12}{24}\selectfont

\CompileMatrices %Compile xymatrix diagrams for faster Latexing

\title[Projective pairs]{Projective pairs of profinite groups}%
\author{Lior Bary-Soroker}%
\address{Institute of Mathematics, Hebrew University, Giv'at-Ram, Jerusalem, 91904, Israel}%
\email{barylior@post.tau.ac.il}%

\thanks{This work is a part of the author's PhD thesis done at Tel Aviv University under the supervision of Prof.\ Dan Haran }%

\date{\today}%
%\dedicatory{}%
%\commby{}%
% ----------------------------------------------------------------
\begin{abstract}
We generalize the notion of a projective profinite group to a
projective pair of a profinite group and a closed subgroup.

We establish the connection with Pseudo Algebraically Closed (PAC)
extensions of PAC fields: Let $M$ be an algebraic extension of a PAC
field $K$. Then $M/K$ is PAC if and only if the corresponding pair
of absolute Galois groups $(\gal(M),\gal(K))$ is projective.
Moreover any projective pair can be realized as absolute Galois
groups of a PAC extension of a PAC field.

Using this characterization we construct new examples of PAC
extensions of relatively small fields, e.g.\ unbounded abelian
extensions of the rational numbers.
\end{abstract}
\maketitle
% ----------------------------------------------------------------
\section{Introduction}
\subsection{Projective pairs}
In the profinite category, or more generally in the pro-$\C$
category for some Melnikov formation of finite groups $\C$ (see
Section~\ref{sec:background}), the projectivity is determined via
$\C$ embedding problems (abbreviated as $\C$-EP). Namely a pro-$\C$
group $\Lambda$ is $\C$-projective if and only if every $\C$-EP,
that is, a diagram
\[
\xymatrix{%
    &\Lambda\ar[d]^{\mu}\\
G\ar[r]^{\alpha}
    &A
}%
\]
in which $\mu, \alpha$ are surjective and $G,A\in \C$ (in particular
$G,A$ are finite), is weakly solvable, i.e., there exists a
homomorphism $\theta \colon \Lambda \to G$ such that $\alpha \theta
= \mu$.

The object that this study concerns is a \emph{$\C$-projective pair}
$(\Gamma, \Lambda)$. Here $\Gamma \leq \Lambda$ are pro-$\C$ groups
with the property that every \emph{$\C$ double embedding problem}
(in short $\C$-DEP) for the pair $(\Gamma, \Lambda)$ is \emph{weakly
solvable}. Roughly speaking a $\C$-DEP is a pair of two $\C$-EPs,
one for $\Gamma$ and one for $\Lambda$, which are compatible. A weak
solution of a $\C$-DEP is a pair of compatible weak solutions of the
corresponding $\C$-EPs. We drop the $\C$ notation, if $\C$ is the
family of all finite groups. See Section~\ref{sec:DEP} for precise
definitions.

The notion of $\C$-projective pairs generalizes the notion of
$\C$-projective groups
(Proposition~\ref{prop:Cproj-simplecharac}). Moreover we give
several characterizations of $\C$-projective pairs including the
\emph{lifting property} (Proposition~\ref{prop:lifting}) and a
non-abelian cohomology characterization
(Proposition~\ref{prop:cohomology}).

\subsection{Projective pairs and pseudo algebraically closed
extensions of fields} The motivation for this new notion of
projective pairs lies in the theory of fields. To explain this
connection we start with the classical case of projective groups.
Ax-Lubotzky-v.d.~Dries Theorem asserts that the class of all
projective groups coincides with the class of all absolute Galois
groups of a special kind of fields, namely Pseudo Algebraically
Closed (PAC) fields, see \cite[Corollary~23.1.3]{FriedJarden2005}.
It is important to note that there are non-PAC fields whose absolute
Galois group is projective, e.g.\ $\bbF_q$,
$\bbC(t)$, and $\bbQ^{\ab}$ (the latter being the maximal abelian
extension of $\bbQ$).

In \cite{JardenRazon1994} Jarden and Razon generalize the notion of
PAC fields and define PAC extensions. (See the introduction of
\cite{Bary-SorokerPACext} for a short survey on PAC extensions and
their applications). Basing on \cite{Bary-SorokerPACext} we prove an
analogous connection between projective pairs and PAC extensions of
PAC fields (see Theorem~\ref{prop:projpair-PACext} below). Note that
in the case $M/K$ is algebraic we have a characterization.

\begin{theorem}\label{prop:projpair-PACext}
\begin{enumerate}
\item
\label{prop:projpair-PACext_1} Let $M$ be a PAC extension of a PAC
field $K$. Then the pair $(\gal(M), \gal(K))$ is projective.
\item
Let $M$ be an algebraic extension of a PAC field $K$. Then $M/K$ is
PAC if and only if the restriction map $(\gal(M),\gal(K))$ is
projective.
\item
Let $(\Gamma, \Lambda)$ be a projective pair. Then there exists a
separable algebraic PAC extension $M$ of a PAC field $K$ such that
$\Gamma\cong \gal(M)$, $\Lambda\cong \gal(K)$.
\end{enumerate}
\end{theorem}

Note that \eqref{prop:projpair-PACext_1} implicitly implies that the
restriction map $\gal(M) \to \gal(K)$ is injective even if $M/K$ is
not algebraic. This is indeed true, see
\cite[Theorem~4.2]{Bary-SorokerPACext}.

In \cite{JardenRazon1994} Jarden and Razon prove that if $K$ is a
countable Hilbertian field and $e\geq 1$ an integer, then for almost
all $\bfsigma = (\sigma_1, \ldots, \sigma_e) \in \gal(K)^e$
\[
K_s(\bfsigma) = \{ x\in K_s \mid \ \forall i\: \sigma_i(x) = x\}
\]
is a PAC extension of $K$. Moreover they prove that if $M/K$ is PAC
and $L/K$ is algebraic, then $LM/L$ is PAC. From these two results
all the known examples of PAC extensions are derived, cf.\
\cite{Bary-SorokerKelmer} for several explicit constructions of that
kind. However much is unknown, for example, for a finitely generated
infinite field $K$ we do not know if there exists a PAC extension
$M/K$ whose absolute Galois group $\gal(M)$ is not finitely
generated \cite[Conjecture~7]{Bary-SorokerJarden}.

We purpose here a new group theoretic method to construct PAC
extensions.
\begin{itemize}
\item Start with a PAC extension $M/K$.
\item Find a PAC extension $E/M$.
Since $M$ is a PAC field, to find $E$ is the purely group theoretic
problem of finding a subgroup $\Gamma$ of $\gal(M)$ such that
$(\Gamma,\gal(M))$ is projective
(Theorem~\ref{prop:projpair-PACext}).
\item By the transitivity of PAC extensions \cite[Theorem~5]{Bary-SorokerPACext} $E/K$ is PAC.
\end{itemize}

Many constructions can be generate by this method. In here we apply
it to relatively small infinite extensions of any countable
Hilbertian field, such as $\bbQ$.

\begin{theorem}
Let $P$ be a countably generated projective group, $K_0$ a countable
Hilbertian field, and $K/K_0$ an abelian extension of $K_0$ such
that
\[
\{ \ord(\sigma) \mid \sigma \in \gal(K/K_0)\} \subseteq \bbN \cup
\{\infty\}
\]
is unbounded. Then there exists a PAC extension $M/K$ such that
$\gal(M) \cong P$.
\end{theorem}

At the moment our method does not apply for a finitely generated
infinite field  $K$. If $(\Gamma,\Lambda)$ is a projective pair,
then $\Gamma$ is a quotient of $\Lambda$
(Corollary~\ref{cor:semidirectproductpp}). Thus if $\gal(M)$ is
finitely generated, then $\gal(E)$ constructed by the above method
will also be finitely generated.

\subsection{The structure of projective pairs}\label{sec:intstruc}
This work also contains some structural study of $\C$-projective
pairs. For example we prove that if $(\Gamma, \Lambda)$ is
$\C$-projective and $\Gamma\neq \Lambda$, then
\begin{enumerate}
\item \label{propzerty1} the normal core of $\Gamma$ is trivial, i.e.\ $\bigcap_{\sigma \in
\Lambda} \Gamma^\sigma = 1$,
\item $(\Lambda:\Gamma) = \infty$, and
\item \label{propzerty3} $\Lambda = N \rtimes \Gamma$, for some normal subgroup
$N$ of $\Lambda$.
\end{enumerate}

%-------------------
It is interesting to note that the analogs (via Galois
correspondence) of some of the properties of projective pairs are
already known for PAC extension. By
Theorem~\ref{prop:projpair-PACext}, the results for PAC fields carry
over to projective pairs. (The opposite implication works only if
$K$ is PAC.)

Nevertheless we bring here group theoretic proofs for several
reasons. First the aesthetic reason -- the group theoretic proofs
are easier. The generality reason -- going via
Theorem~\ref{prop:projpair-PACext} only applies to $\C = $ all
finite groups. Finally the strength reason -- the results about
projective pairs are usually stronger than the corresponding field
theoretic analogs. For example the analog of (c) for a PAC extension
$M/K$ says that if $G$ is a \emph{finite} quotient of $\gal(M)$ that
\emph{regularly occurs} over $K$, then $G$ is a quotient of
$\gal(K)$ \cite[Corollary~6.1]{Bary-SorokerPACext}.

\section{Definition and characterizations of projective pairs}
\label{sec:ppDefinition}
\subsection{Melnikov formations}\label{sec:background}
Throughout this work $\C$ is a fixed Melnikov formation of finite
groups. That means that $\C$ is a family of finite groups that is
closed under taking fiber products and given a short exact sequence
\[
\xymatrix@1{%
1\ar[r]& A \ar[r]& B \ar[r]& C\ar[r]&1
}%
\]
we have that $A,C\in \C$ if and only if $B\in \C$. In particular,
$\C$ is closed under direct products.

The following three families are Melnikov. The family of all finite
groups; the family of all $p$-groups; the family of all solvable
groups. More generally, if $\calS$ is a set of simple finite groups,
then the family of all finite groups whose composition factors are
in $\calS$ is a Melnikov formation.

\subsection{Double embedding problems}\label{sec:DEP}
Let $\Gamma\leq \Lambda$ be pro-$C$ groups. A \textbf{$\C$ double
embedding problem}, or in short \textbf{$\C$-DEP}, for the pair
$(\Gamma, \Lambda)$ is a commutating diagram
\begin{eqnarray}
\label{double embedding problem}%
\xymatrix@C=40pt@R=10pt{%
    &\Lambda\ar[dr]^{\nu}\ar@{.>}[dl]_\theta\\
H\ar[rr]^(.4)\beta
        &&B\\
    &\Gamma\ar[dr]^{\mu}\ar'[u]^{\phi}[uu]\ar@{.>}[dl]_{\eta}\\
G\ar[rr]^{\alpha}\ar[uu]^{j}
        &&A,\ar[uu]^{i}
}%
\end{eqnarray}
where $G,H,A,B\in \C$, $A\leq B$, $G\leq H$, $i,j,\phi$ are the
inclusion maps, and $\alpha,\mu,\beta,\nu$ are surjective. Therefore
a $\C$-DEP consists of two compatible $\C$-EPs: the \textbf{lower}
embedding problem $(\mu,\alpha)$ for $\Gamma$ and the
\textbf{higher} embedding problem $(\nu,\beta)$ for $\Lambda$.

In case $\C$ is the family of all finite groups, we omit the $\C$
notation and simply say that \eqref{double embedding problem} is a
DEP. Sometimes we abbreviate \eqref{double embedding problem} and
write $((\mu,\alpha),(\nu,\beta))$.

A $\C$-DEP is said to be \textbf{split} if $\alpha$ and $\beta$ have
sections,  i.e., there exist $\alpha'\colon A\to G$ and
$\beta'\colon B\to H$ for which $\alpha\alpha' = \id_A$ and
$\beta\beta'=\id_B$. We emphasize that no compatibility condition on
$\alpha'$ and $\beta'$ is required, i.e.\ we allow that $j\alpha'
\neq \beta' i$.

If the groups $G,H,A,B$ are pro-$\C$, then \eqref{double embedding
problem} is a \textbf{pro-$\C$-DEP}.

Given weak solution $\eta\colon \Gamma\to G$ of the lower embedding
problem and weak solution $\theta\colon \Lambda\to H$ of the higher
embedding problem, we say that $(\eta,\theta)$ is a \textbf{weak
solution} of \eqref{double embedding problem} if $\eta$ and $\theta$
are compatible, i.e.\ $\eta = \theta|_{\Gamma}$. Note that
$(\eta,\theta)$ is completely determined by $\theta$: A weak
solution $\theta$ of the higher embedding problem induces a weak
solution of \eqref{double embedding problem} if and only if
$\theta(\Gamma)\leq G$.

\subsection{The definition of $\C$-projective pairs}

\begin{definition}\label{def:pp}
A pair $(\Gamma, \Lambda)$ of pro-$\C$ groups is called
\textbf{$\C$-projective} if any $\C$-DEP is weakly solvable.
\end{definition}

As in the case of $\C$-projective groups
\cite[Lemma~22.3.2]{FriedJarden2005}, the solvability property
extends to pro-$\C$-DEPs.

\begin{proposition}
\label{prop infinite EPs}%
Any pro-$\C$-DEP for a $\C$ projective pair $(\Gamma,\Lambda)$ is
weakly solvable.
\end{proposition}

\begin{proof}
In order to solve pro-$\C$-DEPs for $(\Gamma,\Lambda)$ we need to
solve more general pro-$\C$-DEPs, in which the maps of the
pro-$\C$-DEP are not necessarily surjective.

In case of $\C$-DEPs, we can solve such $\C$-DEPs. Indeed, assume that in \eqref{double embedding problem} $\nu,\mu$ are
not surjective. First $\ker(\alpha),\ker(\beta)\in \C$ since they are normal subgroups. Next $\nu(\Gamma),\mu(\Lambda)\in
\C$ since $\Gamma,\Lambda$ are pro-$\C$. Finally $\alpha^{-1}(\nu(\Gamma)),\beta^{-1}(\mu(\Lambda))\in\C$ follows from
the exact sequences
\[
\xymatrix{%
1\ar[r]&\ker(\alpha)\ar[r]&\alpha^{-1}(\nu(\Gamma))\ar[r]^(0.6)\alpha&\nu(\Gamma)\ar[r]&1\\
1\ar[r]&\ker(\beta)\ar[r]&\beta^{-1}(\mu(\Lambda))\ar[r]^(0.6)\beta&\mu(\Lambda)\ar[r]&1.
}%
\]
Replace $A,B$ with $\nu(\Gamma),\mu(\Lambda)$ and $G,H$ with $\alpha^{-1}(\nu(\Gamma)),\beta^{-1}(\mu(\Lambda))$. In this
new $\C$-DEP all the maps are surjective. So by assumption there is a weak solution.

Let us move to the more general case of pro-$\C$-DEP: Consider a
pro-$\C$-DEP \eqref{double embedding problem} and write $K =
\ker(\beta)$. We prove the assertion in two steps.

\vspace{10pt}%
\noindent\textsc{Step A:} \textit{Finite Kernel.} Assume $K$ is
finite. Then $G$ is open in $KG$ since $(KG:G)\leq |K|$. Choose an
open normal subgroup $U \leq H$ for which $U\cap KG \leq G$ and
$K\cap U = 1$ (note that $K$ is finite and $H$ is Hausdorff). Then
$U\cap KG = U\cap G$. By the second isomorphism theorem (in the
group $UG$) we have
\[
(KG\cap UG:G) = (U\cap (KG\cap UG): U\cap G)=(U\cap KG:U\cap G)=1,
\]
that is to say
\begin{equation}
\label{eq KG intersection UG = G}
   (K G) \cap (U G) = G.
\end{equation}
Write $\Hgag = H/U$, let $\pi \colon H\to \Hgag$ be the quotient
map, $\Ggag = \pi(G)$, $\Bgag = B/\beta(U)$, $\Agag=A/A\cap\beta(U)$
and $\betagag\colon \Hgag\to \Bgag$, $\alphagag\colon \Ggag\to
\Agag$ the epimorphisms induced from $\beta$, $\alpha$,
respectively.

Since $\Hgag\in \C$, there is a homomorphism $\thetagag\colon
\Lambda \to \Hgag$ with $\thetagag(\Gamma)\leq \Ggag$ (let
$\etagag=\thetagag|_{\Gamma}$) for which
\[
\xymatrix{%
    &\Gamma\ar[d]^{\nu}\ar@{.>}[ddl]^(0.7){\etagag}\\
G\ar[r]^{\alpha}\ar[d]
    & A\ar[d]\\
\Ggag\ar[r]^-{\alphagag}
    &\Agag
}%
\qquad
\xymatrix{%
            &&&\Lambda\ar[d]^{\mu}\ar@{.>}[ddl]^(0.7){\thetagag}\\
1\ar[r]
    & K \ar[r]\ar@{=}[d]
        & H \ar[r]^{\beta}\ar[d]_{\pi}
            &B \ar[d] \ar[r]
                &1\\
1\ar[r]
    & K \ar[r]
        & \Hgag \ar[r]^{\betagag}
            &\Bgag \ar[r]
                &1
}%
\]
are commutative diagrams. The right square in the right diagram is a
cartesian square, since $K\cap U = 1$ (\cite[Example
22.2.7(c)]{FriedJarden2005}). Hence we can lift $\thetagag$ to
$\theta \colon \Lambda \to H$ such that $\beta \theta=\mu$
(\cite[Lemma 22.2.1]{FriedJarden2005}). We claim that
$\theta(\Gamma)\leq G$. Indeed,
\[
A \geq \mu(\Gamma) = \beta(\theta(\Gamma)),
\]
hence $\theta(\Gamma) \leq K \beta^{-1}(A) = K \alpha^{-1}(A) = KG$.
Also,
\[
\Ggag \geq \thetagag(\Gamma) = \pi(\theta(\Gamma)),
\]
hence $\theta(\Gamma)\leq UG$. Then, from \eqref{eq KG intersection
UG = G} we have $\theta(\Gamma) \leq (KG) \cap (UG) = G$, as
claimed.

\vspace{10pt}%
\noindent\textsc{Step B:} \textit{The General Case.} We use Zorn's
Lemma. Consider the family of pairs $(L,\theta)$ where $L\subseteq
K$ is normal in $H$, $\theta$ is a weak solution of the following
embedding problem, and $\theta(\Gamma) \subseteq GL/L$.
\[
\xymatrix{ &&&\Lambda\ar[d]\ar[ld]_{\theta}
\\
1\ar[r]&K/L\ar[r]& H/L\ar[r]^-{\bar\beta}&B\ar[r]&1. }
\]
We say that $(L,\theta) \leq (L',\theta')$ if $L\subseteq L'$ and
\[
\xymatrix{ &\Lambda\ar[d]\ar[ld]_{\theta}\ar[ldd]^(.6){\theta'}
\\
H/L\ar[r]|{~}\ar[d]& B\ar@{=}[d]\\
H/L'\ar[r]& B}
\]
is commutative. For a chain $\{(L_i,\theta_i)\}$ we define a lower
bound $(L,\theta)$ by $L = \bigcap_{i}L_i$ and $\theta = \invlim
\theta_i$ (note that $\theta(\Gamma)\subset GL/L$ by
\cite[Lemma~1.2.2(b)]{FriedJarden2005}). By Zorn's Lemma there
exists a minimal element $(L,\theta)$ in the family. We claim that
$L = 1$. Otherwise, there is an open normal subgroup $U$ of $H$ with
$L \not\leq U$. Part A gives (since $L/U\cap L$ is finite) a weak
solution $\theta'$ of the following embedding problem such that
$\theta'(\Gamma) \subseteq G(U\cap L)/(U\cap L)$.
\[
\xymatrix{ &&&\Lambda\ar[d]_{\theta}\ar[dl]_{\theta'}
\\
1\ar[r]&L/U\cap L\ar[r] &H/U\cap L \ar[r]& H/L. }
\]
Hence $(L,\theta)$ is not minimal.
\end{proof}

\begin{remark}
For an algebraic PAC extension $M/K$, special kind of finite double embedding problems for $(\gal(M),\gal(K))$ are weakly
solvable \cite[Proposition~4.5]{Bary-SorokerPACext}. However it is unknown whether one can solve profinite double
embedding problems for $(\gal(M),\gal(K))$. Hence Proposition~\ref{prop infinite EPs} strengthens this property in case
$K$ is PAC (via Theorem~\ref{prop:projpair-PACext}).
\end{remark}

The next result shows that $\C$-projective pairs generalize
$\C$-projective groups.

\begin{proposition}\label{prop:Cproj-simplecharac}
A pro-$\C$ group $\Lambda$ is $\C$-projective if and only if the
pair $(1, \Lambda)$ is $\C$-projective.
\end{proposition}

\begin{proof}
A $\C$-EP $(\mu\colon \Lambda \to A, \alpha\colon G\to A)$ for
$\Lambda$ is weakly solvable if and only if the $\C$-DEP $((1\to 1,
1\to 1), (\mu,\alpha))$ for $(1,\Lambda)$ is weakly solvable. This
implies the first equivalence.

Let $\beta\colon H\to A$ be any epimorphism of pro-$\C$-group
satisfying $G\leq H$ and $\beta|_G = \alpha$. Then $\theta$ is a
weak solution of $(\mu,\alpha)$ if and only if $(\theta,\theta)$ is
a weak solution of $((\mu,\alpha), (\mu,\beta\colon H \to A))$.
\end{proof}

\begin{remark}
From Proposition~\ref{prop:transitivity} below it also follows that
$\Lambda$ is $\C$-projective if and only if $(\Lambda,\Lambda)$ is
$\C$-projective.
\end{remark}

From technical perspective, it is important to dominate a
pro-$\C$-DEP by a more convenient one, e.g.\ split pro-$\C$-DEP. Let
us make a precise definition.

\begin{definition}
Let $\Gamma\leq \Lambda$ be pro-$\C$ groups and consider the
following two $\C$-DEP for $(\Gamma,\Lambda)$.
\[
\xymatrix@C=40pt@R=10pt{%
    &\Lambda\ar[dr]^{\nu}
                &&&\Lambda\ar[dr]^{\nuhat}\\
H\ar[rr]^(.4)\beta
        &&B
            &\Hhat\ar[rr]^(.4)\betahat
                    &&\Bhat\\
    &\Gamma\ar[dr]^{\mu}\ar'[u]^{\phi}[uu]
                &&&\Gamma\ar[dr]^{\muhat}\ar'[u]^{\phi}[uu]\\
G\ar[rr]^{\alpha}\ar[uu]^{j}
        &&A\ar[uu]^{i}
            &\Ghat\ar[rr]^{\alphahat}\ar[uu]^{\jhat}
                    &&\Ahat.\ar[uu]^{\ihat}
}%
\]
We say that $((\muhat,\alphahat),(\nuhat,\betahat))$ \textbf{dominates} $((\mu,\alpha),(\nu,\beta))$ if there exist
epimorphisms $\pi_i$, $i=1,2,3,4$ making the following diagram commutate.
\[
\xymatrix{%
    &\Hhat\ar[rr]^{\betahat}\ar'[d][dd]^{\pi_2}
           &&\Bhat\ar[dd]^{\pi_4}\\
\Ghat\ar[rr]^(.6){\alphahat}\ar[dd]^{\pi_1}\ar[ur]^{\jhat}
        &&\Ahat\ar[dd]_(.7){\pi_3}\ar[ur]^{\ihat}\\
    &H\ar'[r][rr]^\beta
        &&B\\
G\ar[rr]^\alpha \ar[ur]^j
        &&A.\ar[ur]^i
}%
\]
Clearly every weak solution $(\etahat,\thetahat)$ of the dominating
$\C$-DEP induces a solution $(\eta,\theta)$ of the dominated
$\C$-DEP by setting $\eta = \pi_1\etahat$ and $\theta =
\pi_2\thetahat$.
\end{definition}

\begin{lemma}\label{lem:DEPdominated}
Consider a $\C$-DEP \eqref{double embedding problem} for a pair
$(\Gamma, \Lambda)$ of pro-$\C$ groups. Assume that both the higher
and lower embedding problems are weakly solvable. Then \eqref{double
embedding problem} is dominated by a split $\C$-DEP.
\end{lemma}

\begin{proof}
Let $\theta\colon \Lambda \to H$ be a weak solution of the higher
embedding problem and $\eta\colon \Gamma\to G$ a weak solution of
the lower embedding problem. Choose an open normal subgroup $N\leq
\Lambda$ such that $N\leq \ker(\theta)$ and $\Gamma\cap N\leq
\ker(\eta)$.

Let $\Bhat=\Lambda/N$, $\Ahat=\Gamma/\Gamma\cap N$ and let $\Hhat = H\times_B \Bhat$, $\Ghat = G\times_A \Ahat$. Then the
upper rows in the following commutative diagrams define a dominating $\C$-DEP.
\[
\xymatrix{%
    &\Gamma\ar[d]_{\muhat}\ar@/^/[dd]^{\mu}\\
\Ghat\ar[r]^{\alphahat}\ar[d]^{\pi_1}
    &\Ahat\ar[d]\\
G\ar[r]^\alpha
    &A
}%
\qquad
\xymatrix{%
    &\Lambda\ar[d]_{\nuhat}\ar@/^/[dd]^{\nu}\\
\Hhat\ar[r]^{\betahat}\ar[d]^{\pi_2}
    &\Bhat\ar[d]\\
H\ar[r]^\beta
    &B
}%
\]
(Here all the maps are canonically defined.)

To finish the proof we need to show that both $\alphahat$ and
$\betahat$ defined in the above diagram have sections. Let
$\alphahat' \colon \Ahat\to \Ghat$ be defined by $\alphahat'(x) =
(\eta(x),x)$, $x\in \Ahat$, and similarly let $\betahat'\colon \Bhat
\to \Hhat$ be defined by $\betahat'(x) = (\theta(x),x)$, $x\in
\Bhat$. Then, $\alphahat (\alphahat'(x)) = x$, $x\in \Ahat$ and
$\betahat (\betahat'(x)) = x$, $x\in \Bhat$, i.e.\ both $\alphahat$
and $\betahat$ split, as needed.
\end{proof}

\begin{corollary}\label{cor:ppreductiontodoublesplit}
Let $(\Gamma,\Lambda)$ be a pair of pro-$\C$ groups and suppose that
$\Lambda$ is $\C$-projective. Then $(\Gamma,\Lambda)$ is
$\C$-projective if and only if every split $\C$-DEP is weakly
solvable.
\end{corollary}

\begin{proof}
Since $\Lambda$ is $\C$-projective, $\Gamma$ is also
$\C$-projective. In other words, every finite embedding problem for
$\Lambda$ (resp.\ $\Gamma$) is weakly solvable.
Lemma~\ref{lem:DEPdominated} implies that every $\C$-DEP for
$(\Gamma,\Lambda)$ is dominated by a split $\C$-DEP.

The converse is trivial.
\end{proof}

Recall that for a pro-$\C$ group $\Lambda$ to be $\C$-projective it is necessary and sufficient that any short exact
sequence of pro-$\C$ groups
\[
\xymatrix{1\ar[r]&K\ar[r]& \Delta\ar[r]& \Lambda\ar[r]&1}
\]
splits. A similar characterization is given in the next result for a pair $(\Gamma,\Lambda)$ of pro-$\C$ groups.

\begin{corollary}\label{cor:splitting}
Let $(\Gamma,\Lambda)$ be a pair of pro-$\C$ groups. Then
$(\Gamma,\Lambda)$ is $\C$-projective if and only if the rows of any
exact commutative diagram of pro-$\C$ groups
\[
\xymatrix{%
\Delta \ar[r]^\beta
    & \Lambda \ar[r]\ar@<1ex>@{.>}[l]^{\beta'}
        &1 \\
E \ar[r]^\alpha \ar[u]^{\psi}
    & \Gamma\ar[r]\ar[u]_{\phi}\ar@<1ex>@{.>}[l]^{\alpha'}
        &1
        \\
1\ar[u]&1\ar[u]
}%
\]
\emph{compatibly} split. That is to say, there exists a section
$\beta'\colon \Lambda\to \Delta$ of $\beta$ such that $\beta' \phi
(\Gamma)\leq \psi(E)$ and $\alpha'$ defined by
$\psi\alpha'=\beta'\phi$ is a section of $\alpha$.
\end{corollary}

\begin{proof}
Since the sections of an epimorphism $\gamma\colon M\to N$
correspond bijectively to solutions of the embedding problem
$(\id\colon N\to N, \gamma\colon M\to N)$, the assertion follows
immediately from Proposition~\ref{prop infinite EPs}.
\end{proof}

\subsection{The lifting property}

The following lifting property is a key property in proving the structural results mentioned in
Section~\ref{sec:intstruc}. This is a stronger version of the lifting property for PAC extensions, since it applies to
pro-$\C$-DEPs, and not only to finite DEPs.

\begin{proposition}[The lifting property]\label{prop:lifting}
Let $(\Gamma,\Lambda)$ be $\C$-projective and consider a pro-$\C$-DEP
\begin{equation} \label{eq:depinsideproof}
\xymatrix@R=10pt{%
    &\Lambda\ar[dr]^{\nu}\\
H\ar[rr]^(.4)\beta
        &&B\\
    &\Gamma\ar[dr]^{\mu}\ar'[u][uu]\ar@{.>}[dl]_{\eta}\\
G\ar[rr]^{\alpha}\ar[uu]
        &&A,\ar[uu]
}%
\end{equation}
for $(\Gamma,\Lambda)$. Then any weak solution $\eta$ of the lower embedding problem can be lifted to a weak solution
$(\eta,\theta)$ of \eqref{eq:depinsideproof}.
\end{proposition}

\begin{proof}

Let $\eta\colon \Gamma\to G$ be a weak solution of the lower
embedding problem $(\mu,\alpha)$. Define $\Hhat = H \times_B
\Lambda$ and let $\pi\colon \Hhat \to H$ and $\betahat\colon
\Hhat\to \Lambda$ be the quotient maps. Let $\Ghat
=\{(\eta(\gamma),\gamma )\mid \gamma\in \Gamma\}\leq G\times_A
\Gamma$ and $\alphahat=\betahat|_{\Ghat}$. Then
$\alphahat((\eta(\gamma),\gamma))=\gamma$, for all $\gamma\in
\Gamma$, and hence $\alphahat$ is an isomorphism. Thus
$\alphahat^{-1}$ is the \emph{unique} weak solution of the lower
embedding problem of
\[
\xymatrix@R=10pt{%
    & \Lambda\ar[rd]^{\id}\\
\Hhat\ar[rr]^(.4){\betahat}
        && \Lambda\\
    & \Gamma\ar'[u][uu]\ar[dr]^{\id}\ar[dl]_(.3){\alphahat^{-1}}\\
\Ghat\ar[rr]^{\alphahat}\ar[uu]
        && \Gamma.\ar[uu]
}%
\]
By Proposition~\ref{prop infinite EPs} there exists a weak solution
$(\etahat,\thetahat)$ of the above pro-$\C$-DEP. Hence, $\etahat =
\alphahat^{-1}$.

Let $\eta' = \pi|_{\Ghat}\etahat$ and $\theta = \pi \thetahat$. Then
$(\eta',\theta)$ is a weak solution of \eqref{eq:depinsideproof}.
Moreover \[\eta'(\gamma) = \pi(\alphahat^{-1}(\gamma)) = \pi
((\eta(\gamma),\gamma )) = \eta(\gamma),\] i.e.\ $(\eta,\theta)$ is
a weak solution of \eqref{eq:depinsideproof}, as needed.
\end{proof}

We give two characterizations of $\C$-projective pairs. The first
follows from the lifting property using the same argument that
implied Corollary~\ref{cor:splitting} from Proposition~\ref{prop
infinite EPs}. The second is in terms of non-abelian cohomology.

\begin{corollary}\label{cor:lifting_splitting}
Let $(\Gamma,\Lambda)$ be a $\C$-projective pair and consider a
diagram as in Corollary~\ref{cor:splitting}. Then any section
$\alpha'$ of $\alpha$ can be lifted to a section $\beta'$ of
$\beta$.
\end{corollary}

\begin{proposition}\label{prop:cohomology}
Let $\Lambda$ be a $\C$-projective group and $\Gamma$ a subgroup.
The pair $(\Gamma, \Lambda)$ is $\C$-projective if and only if for
any pro-$\C$ group $A$ on which $\Lambda$ acts
%$H^2(\Lambda,A)=0$ and
the restriction map
\[
H^1(\Lambda, A) \to H^1(\Gamma, A)
\]
is surjective.
\end{proposition}

\begin{proof}
Recall that there is a natural identification between $H^1(\Lambda,
A)$ and sections of the quotient map $\beta\colon A\rtimes \Lambda
\to \Lambda$. More precisely, every $x\in H^1(\Lambda, A)$ induces
the section $\beta'$ defined by $\beta'(\lambda) = x(\lambda)
\lambda$, $\lambda\in \Lambda$ and vice versa a section $\beta'$ of
$\beta$ induces $x\in H^1(\Lambda, A)$ defined by $x (\lambda) =
\beta'(\lambda)\lambda^{-1}$.

Assume $(\Gamma, \Lambda)$ is $\C$-projective. Let $x\in H^1(\Gamma,
A)$. It defines an embedding
\[
\psi \colon \Gamma\to A\rtimes \Gamma \leq A\rtimes \Lambda, \qquad
\psi (\gamma) = x(\gamma) \gamma,\ \forall \gamma\in \Gamma.
\]
Corollary~\ref{cor:lifting_splitting} applied to
\[
\xymatrix{%
A\rtimes \Lambda \ar[r]^{\beta}&
    \Lambda\ar[r]&1\\
\Gamma\ar[r]^{\id} \ar[u]^{\psi}&
    \Gamma\ar[u]\ar[r]&1
}%
\]
(where $\beta$ is the quotient map) gives $\beta'\colon \Lambda \to
A\rtimes \Lambda$ such that $\beta \beta'$ is the identity on
$\Lambda$ and $\beta'|_{\Gamma} = \psi$. Then $\beta'$ induces $y\in
H^1(\Lambda,A)$ defined by $y(\lambda) = \beta'(\lambda)
\lambda^{-1}$. Clearly the restriction of $y$ to $\Gamma$ is $x$.

Next assume that the restriction map $H^1(\Lambda, A) \to
H^1(\Gamma, A)$ is surjective. Consider a commutative exact diagram
\[
\xymatrix{%
1\ar[r]&B\ar[r]&\Delta \ar[r]^\beta
    & \Lambda \ar[r]%\ar@<1ex>@{.>}[l]^{\beta'}
        &1 \\
1\ar[r]&A\ar[r]\ar[u]&E \ar[r]^\alpha \ar[u]^{\psi}
    & \Gamma\ar[r]\ar[u]_{\phi}%\ar@<1ex>@{.>}[l]^{\alpha'}
        &1\ar@{}[r]_(0.1).&
        \\
&&1\ar[u]&1\ar[u]
}%
\]
Since $\Lambda$ is $\C$-projective, $\Gamma$ is also
$\C$-projective, and hence both rows split. Identify $\Delta$ with
$B\rtimes \Lambda$ via some fixed section of $\beta$. Let $\alpha'$
be a section of $\alpha$ and $x\in H^1(\Gamma,A) \leq H^1(\Gamma,B)$
the corresponding cocycle (i.e.\ $x(\gamma) =
\alpha'(\gamma)\gamma^{-1}$). By assumption there exists a cocyle
$y\in H^1(\Lambda, B)$ satsifying $y|_\Gamma = x$. Let $\beta'$ be
the induced section of $\beta$. Then for all $\gamma\in \Gamma$ we
have
\[
\beta'(\gamma) = y(\gamma) \gamma = x(\gamma) \gamma =
\alpha'(\gamma).
\]
Therefore $\beta'|_\Gamma = \alpha'$ and by
Corollary~\ref{cor:lifting_splitting} $(\Gamma,\Lambda)$ is
$\C$-projective.
\end{proof}

\begin{proposition}[Transitivity]\label{prop:transitivity}
Let $\Lambda_3\leq \Lambda_2\leq \Lambda_1$ be pro-$\C$ groups.
\begin{enumerate}
\item
    If $(\Lambda_3,\Lambda_1)$ is $\C$-projective, then $(\Lambda_3,\Lambda_2)$ is
    $\C$-projective.
\item
    If $(\Lambda_3,\Lambda_2)$ and $(\Lambda_2,\Lambda_1)$ are $\C$-projective, then so is
    $(\Lambda_3,\Lambda_1)$.
\end{enumerate}
\end{proposition}

\begin{proof}
For (a) assume that $(\Lambda_3,\Lambda_1)$ is $\C$-projective. Consider a commutative diagram
\[
\xymatrix{%
\Lambda_1\ast \Delta_2 \ar[r]^{\alpha_1} \ar@<.5ex>[d]^{\pi}
    &\Lambda_1 \ar[r]
        &1\\
\Delta_2 \ar[r]^{\alpha_2}\ar@<.5ex>[u]^{\psi_2}
    &\Lambda_2 \ar[r]\ar[u]
        &1\\
\Delta_3 \ar[r]^{\alpha_3}\ar[u]^{\psi_3}
    &\Lambda_3\ar[u]\ar[r]
        &1.
}%
\]
Here $\psi_3$ is injective, $\psi_2$ is the inclusion map,
$\Lambda_1 \ast \Delta_2$ is the free product of $\Lambda_1$ and
$\Delta_2$ in the pro-$\C$ category\footnote{The free product exists
in the pro-$\C$ category. Indeed it is the maximal pro-$\C$ quotient
of the profinite free product of $\Lambda_1$ and $\Lambda_2$},
$\alpha_1$ is defined by $\alpha_1|_{\Lambda_1}= \id$,
$\alpha_1|_{\Delta_2}=\alpha_2$, and $\pi$ is defined by
$\pi|_{\Lambda_1}=1$, $\pi|_{\Delta_2}=\id$.

There exist compatible sections $\beta_3, \beta_1$ of $\alpha_3,
\alpha_1$ (Corollary~\ref{cor:splitting}). Let $\beta_2 = \pi
\beta_1|_{\Lambda_2}$. By the above commutative diagram,
$\beta_3,\beta_2$ are compatible sections of $\alpha_3,\alpha_2$,
and thus $(\Lambda_3,\Lambda_2)$ is $\C$-projective (again
Corollary~\ref{cor:splitting}).

(b) easily follows from Proposition~\ref{prop:cohomology}: Let $A$
be a pro-$\C$ group together with a $\Lambda_1$-action. Since the
restriction map $r_{1,3}\colon H^1(\Lambda_1,A)\to H^1(\Lambda_3,A)$
factors as
\[
\xymatrix{%
H^1(\Lambda_1,A)\ar[r]_{r_{1,2}}\ar@/^10pt/[rr]^{r_{1,3}}
    &H^1(\Lambda_2,A)\ar[r]_{r_{2,3}}
        & H^1(\Lambda_3,A)
}%
\]
and both $r_{1,2}$ and $r_{2,3}$ are surjective (Proposition~\ref{prop:cohomology}) we get that $r_{1,3}$ is surjective.
Consequently, $(\Lambda_3, \Lambda_1)$ is $\C$-projective (again Proposition~\ref{prop:cohomology}).
\end{proof}

\begin{remark}
Let $\Lambda_3\leq \Lambda_2\leq \Lambda_1$ be pro-$\C$ groups. We
show it does not suffice that $(\Lambda_3,\Lambda_1)$ be
$\C$-projective for $(\Lambda_2,\Lambda_1)$ to be $\C$-projective.
For this purpose look at $1 \leq p \bbZ_p \leq \bbZ_p$.

Then $(1, \bbZ_p)$ is $\C$-projective
(Proposition~\ref{prop:Cproj-simplecharac}) while $(p\bbZ_p,\bbZ_p)$
is not (Proposition~\ref{prop:normalprojectivepairs} below).
\end{remark}

Projective pairs behave well under taking subgroups.
\begin{proposition}\label{prop:ppundersubgroups}
Let $( \Gamma, \Lambda)$ be a $\C$-projective pair, let
$\Lambda_0\leq \Lambda$ be a subgroup, and write $\Gamma_0 =
\Lambda_0\cap \Gamma$. Then $( \Gamma_0, \Lambda_0)$ is
$\C$-projective.
\end{proposition}

\begin{proof}
Let $E_0\leq \Delta_0$ be pro-$C$ groups with $\psi\colon E_0 \to
\Delta_0$ the inclusion map, and let $\alpha_0\colon E_0\to
\Gamma_0$ and $\beta_0\colon \Delta_0 \to \Lambda_0$ be epimorphisms
satisfying $\beta_0|_{E_0} = \alpha_0$. Let $\Delta = \Lambda\ast
\Delta_0$, $E = \Gamma_0 \ast E_0$ and let $i_1\colon \Delta_0\to
\Delta$ and $i_2\colon E_0\to E$ be the corresponding injections. We
define maps $\pi_1\colon \Delta \to \Delta_0$ and $\beta\colon
\Delta \to \Delta_0$ by setting $\pi_1|_{\Lambda} = 1$,
$\pi_1|_{\Delta_0} = \id$, $\beta|_{\Gamma} = \id$, and
$\beta|_{\Delta_0} = \beta_0$. Similarly we define $\pi_2\colon E
\to E_0$ and $\alpha\colon E\to \Gamma$.

\[
\xymatrix{%
    &\Delta\ar[rr]^{\beta}\ar@<1ex>[dl]^{\pi_1}
           &&\Lambda\\
\Delta_0\ar[rr]^(.7){\beta_0}\ar[ur]^{i_1}
        &&\Lambda_0\ar[ur]\\
    &E\ar'[r]^(.6)\alpha[rr]\ar'[u]^{\psi}[uu]\ar@<1ex>[dl]^{\pi_2}
        &&\Gamma\ar[uu]\\
E_0\ar[rr]^{\alpha_0} \ar[ur]^{i_2}\ar[uu]^{\psi_0}
        &&\Gamma_0.\ar[ur]\ar[uu]
}%
\]

By Corollary~\ref{cor:splitting} there exist compatible sections
$\alpha'$ and $\beta'$ of $\alpha$ and $\beta$, respectively. Let
$\alpha'_0 = \pi_2 \alpha'$ and $\beta_0' = \pi_1 \beta'$. Then the
above commutative diagram implies that $\alpha'_0$ and $\beta'_0$
are compatible sections of $\alpha_0$ and $\beta_0$, respectively.
Hence $(\Gamma_0,\Lambda_0)$ is $\C$-projective (by the same
corollary).
\end{proof}

\begin{corollary}\label{cor:splittingpp}
Let $(\Gamma, \Lambda)$ be a $\C$-projective pair and let $N\leq
\Gamma$. Then there exists $M\leq \Lambda$ such that $N = \Gamma
\cap M$ and $\Gamma M = \Lambda$. Moreover, if $N\normal \Gamma$,
then $M\normal \Lambda$.
\end{corollary}

\begin{proof}
Let $\Nhat = \bigcap_{\sigma} N^\sigma$ and let $\eta\colon \Gamma
\to \Gamma/\Nhat$ be the natural quotient map. Lift $\eta$ to a
solution $(\eta,\theta)$ of the DEP
\[
((\Gamma\to 1, \Gamma/\Nhat\to 1), (\Lambda\to 1, \Gamma/\Nhat\to
1)).
\]
Let $\Mhat = \ker(\theta)$ and $M = \theta^{-1}(N/\Nhat)$. Then
(since $\eta = \theta|_{\Gamma}$)
\[
N = \eta^{-1}(N/\Nhat) = \Gamma \cap \theta^{-1}(N/\Nhat) =
\Gamma\cap M.
\]
Since $\theta(\Gamma) = \Gamma/\Nhat$, it follows  that $\Gamma
\Mhat = \Lambda$, and in particular, $\Gamma  M = \Lambda$.

To conclude the proof, note that if $N\normal\Gamma$, then
$N=\Nhat$, and hence $M=\Mhat$. So $M\normal \Lambda$, as needed.
\end{proof}

Taking $N=1$ in the above result we get the following splitting
corollary.

\begin{corollary}\label{cor:semidirectproductpp}
If $(\Gamma, \Lambda)$ is $\C$-projective, then $\Lambda \cong M
\rtimes \Gamma$ for some $M\normal \Lambda$.
\end{corollary}

\section{Families of projective pairs}
\subsection{Free products}
We say that $\Gamma$ is a free factor in $\Lambda$ if there exists a
subgroup $N$ of $\Lambda$ such that $\Lambda = \Gamma\ast N$.

\begin{proposition}\label{prop_freefactor}
Let $\Gamma$ be a free factor of a $\C$-projective group $\Lambda$.
Then $(\Gamma, \Lambda)$ is $\C$-projective.
\end{proposition}

\begin{proof}
By assumption $\Lambda = \Gamma\ast N$. Consider a diagram as in
Corollary~\ref{cor:splitting}. Let $M = \beta^{-1} (N) \leq \Delta$
and let $\gamma=\beta|_{N}$. Since $\Lambda$ and $N$ are
$\C$-projective, there exist sections $\alpha',\gamma'$ of
$\alpha,\gamma$, respectively. Let $\beta' \colon \Lambda\to \Delta$
be the induced map. Then $\beta'$ is a section of $\beta$ which is
compatible with $\alpha'$. Thus $(\Gamma,\Lambda)$ is
$\C$-projective (Corollary~\ref{cor:splitting}).
\end{proof}

For a cardinal $\kappa$ let $\Fhat_\kappa$ denote the free pro-$\C$
group. The following result appears in \cite{HaranLubotzky1985} (for
$\kappa = \aleph_0$).

\begin{lemma}[Haran and Lubotzky]\label{lem:ppHL}
Let $\kappa$ be an infinite cardinal and let $P$ be a
$\C$-projective profinite group of rank $\leq \kappa$. Then
$\Fhat_\kappa \cong P\ast \Fhat_{\kappa}$.
\end{lemma}

Combining the above two results yields a family of $\C$-projective
pairs.

\begin{corollary}\label{cor:ppHL}
Let $\kappa$ be an infinite cardinal and $\Lambda = \Fhat_\kappa$.
Then any $\C$-projective group $\Gamma$ of rank $\leq \kappa$ can be
embedded in $\Lambda$ such that $(\Gamma,\Lambda)$ is
$\C$-projective.
\end{corollary}

\subsection{Random finitely generated subgroups}
Let us start by fixing some notation. We write $e$-tuples in bold
face letters, e.g., $\bfb=(b_1, \ldots, b_e)$. For a homomorphism of
profinite groups $\beta\colon H\to B$, we write that
$\beta(\bfh)=\bfb$ if $\beta(h_i) = b_i$ for all $i=1, \ldots, e$.
Let $C$ be a coset of a subgroup $N^e$ in $H^e$, where $N\leq
\ker(\beta)$. By abuse of notation we write $\beta(C)$ for the
unique element $\bfb\in B^e$ such that $\beta(\bfc)=\bfb$ for every
$\bfc \in C$.

For the result of this section we need to add one condition to $\C$,
namely we require that $\C$ be closed under taking subgroups. Then
the pro-$\C$ category is closed under taking subgroups.

\begin{proposition}
Let $\Lambda = \Fhat_\omega$ be the free pro-$\C$ group of countable
rank. Then for almost all $\bfsigma \in \Lambda^e$ (w.r.t.\ the Haar
measure on $\Lambda$) $(\left<\bfsigma\right>, \Lambda)$ is
$\C$-projective.
\end{proposition}

\begin{proof}
Let $m$ denote the normalized Haar measure on $\Lambda^e$. Let
\begin{eqnarray}
\label{eq_ep}%
\xymatrix{%
    &\Lambda\ar[d]^\mu\\
H\ar[r]^\beta
    &B
}%
\end{eqnarray}
be a $\C$-EP for $\Lambda$, let $\bfb \in B^e$, let $A =
\left<\bfb\right>$ be the subgroup of $B$ generated by $\bfb$, and
let $\bfh \in H^e$ be such that $\beta(\bfh) = \bfb$. Define $\Sigma
= \Sigma(\bfb,\bfh,\mu,\beta) \subseteq \Lambda^e$ to be the
following set.
\[
\Sigma=\{\bfsigma \in \Lambda^e\mid (\mu(\bfsigma) = \bfb)\
\Rightarrow\ (\exists \theta\colon \Lambda \to H, \
(\beta\theta=\mu)\wedge (\theta(\bfsigma)=\bfh))\},
\]
that is to say, all $\bfsigma \in \Lambda^e$ such that there exists
a weak solution $\theta$ of \eqref{eq_ep} with $\theta(\bfsigma) =
\bfh$, provided $\mu(\bfsigma) = \bfb$. Note that $\Sigma =
(\Sigma\cap C)\cup (\Lambda^e \smallsetminus C)$, where $C$ is the
coset of $\ker(\mu)^e$ in $\Lambda^e$ for which $\mu(C)=\bfb$.

We break the proof into three parts. In the first two we show that
$m(\Sigma\cap C) = m(C)$, and hence $m(\Sigma) = 1$.

\textsc{Part A:}\emph{Construction of solutions. } Let
\[
\Delta = \{ (b_i)\in H^{\bbN} \mid \beta(b_i) = \beta(b_j)\ \forall
i,j\in \bbN\}.
\]
It is equipped with canonical projections $\pi_i\colon \Delta\to
H_i$, $i\in \bbN$. Set $\betahat\colon \Delta \to B$ by $\betahat(x)
= \beta\pi_i(x)$, $x\in \Delta$. Note that $\betahat$ does not
depend on $i$ and is an epimorphism.

Let $\theta\colon \Lambda \to \Delta$ be a solution of $(\mu\colon
\Lambda \to B, \betahat\colon \Delta\to B)$ (for the existence of
$\theta$ see \cite[Proposition~25.6.2]{FriedJarden2005}). Then for
every $i\in \bbN$ the map $\theta_i = \pi_i\theta$ is a solution of
\eqref{eq_ep}. Moreover, by
\cite[Lemma~2.5]{Bary-SorokerHaranHarbater} the set
$\{\ker(\theta_i)\}$ is an independent set of subgroups of
$\ker(\mu)$.

\vspace{10pt} \textsc{Part B:}\emph{Calculating $m(\Sigma)$. } For
each $i\in \bbN$ take the coset $X_i$ of $\ker(\theta_i)^e$ with
$\theta_i(X_i) = \bfh$. Then, since
\[
\mu(X_i) = \beta \theta_i(X_i) = \beta(\bfh) = \bfb,
\]
it follows that $X_i \subseteq C$. Moreover, Part A implies that
$\{X_i\mid i\in \bbN\}$ is an independent set in $C$.

By the Borel-Cantelli Lemma %\cite[Lemma 18.3.5]{FriedJarden2005},
since $\sum_i m_{C}(X_i) = \sum_i \frac{|B|^e}{|H|^e} = \infty$ we
get that $m_C(X) = 1$. Here $m_C$ is the normalized Haar measure on
$C$ and $X = \cap_{j = 1}^\infty \cup_{i=j}^\infty X_i$. So it
suffices to show that $X\subseteq \Sigma$.

Indeed, let $\bfsigma \in X$. Then $\bfsigma \in X_i$ for some $i$.
It implies that $\theta_i$ is a solution of \eqref{eq_ep} and that
$\theta_i(\bfsigma) =\bfh$. Hence $\bfsigma\in \Sigma$ and
$X\subseteq \Sigma$, as desired.
\vspace{10pt}\textsc{Part C:}\emph{Conclusion. } Let $\Upsilon$ be
the intersection of all $\Sigma(\bfb,\bfh,\mu,\beta)$. Since there
are only countably many of them and each is of measure $1$, we have
$m(\Upsilon)=1$. Let $\bfsigma \in \Upsilon$ and let $\Gamma =
\left<\bfsigma\right>$.

Then $(\Gamma,\Lambda)$ is $\C$-projective. Indeed, consider a
$\C$-DEP as in \eqref{double embedding problem} and choose $\bfh\in
G$ such that $\beta(\bfh)=\mu(\bfsigma)$. Then, since $\bfsigma \in
\Sigma(\mu(\bfsigma), \bfh ,\mu,\beta)$, there exists a homomorphism
$\theta\colon \Lambda \to H$ such that
$\theta(\Gamma)=\left<\theta(\bfsigma)\right>=\left<\bfh\right>\leq
G$.
\end{proof}

\begin{remark}
In the above theorem we actually prove  that for almost all
$\bfsigma\in \Lambda^e$ the pair $(\left<\bfsigma\right> ,\Lambda)$
has the following stronger lifting property. For any $\C$-EP
\eqref{double embedding problem} and for any $\bfh\in G^e$ that
satisfies $\alpha(\bfh)=\mu(\bfsigma)$ there exists a weak solution
$\theta\colon \Lambda\to B$ with $\theta(\bfsigma) = \bfh$.
\end{remark}

\section{Restrictions on projective pairs}

\begin{lemma}\label{lem:normalprojectivepairs}
Let $(\Gamma, \Lambda)$ be a $\C$-projective pair and assume that
$\Gamma\normal \Lambda$. Then either $\Gamma=1$ or $\Gamma=\Lambda$.
\end{lemma}

\begin{proof}
Assume that both $\Gamma$ and $\Lambda/\Gamma$ are not trivial, and
let $\eta\colon \Gamma\to A$ and $\nu \colon \Lambda \to G$ be
epimorphism onto nontrivial $\C$-groups. Recall that the wreath
product of $A$ and $G$, denoted by $A\wr G$, is the semidirect
product $A^G \rtimes G$ w.r.t.\ the translation action of $G$ on
$A^G$. The exact sequence
\[
\xymatrix@1{%
1\ar[r] & A^G \ar[r] & A\wr G\ar[r]^-{\alpha} & G \ar[r]&1,
}%
\]
where $\alpha$ is the quotient map, implies that $A\wr G\in \C$.
Identify $A$ with the subgroup $A^1$ of $A\wr G$.

By the lifting property (Proposition~\ref{prop:lifting}) we can
extend $\eta$ to a weak solution $(\eta,\theta)$ of the $\C$-DEP
\[
\xymatrix@R=10pt{%
    &\Lambda\ar[dr]^{\nu}\ar[dl]_{\theta}\\
A\wr G\ar[rr]^\alpha
        &&G\\
    &\Gamma\ar[dr]\ar'[u][uu]\ar[dl]_{\eta}\\
A\ar[rr]\ar[uu]
        &&1\ar[uu]
}%
\]
Since $\Gamma\normal \Lambda$ we get that $A = \eta(\Gamma) =
\theta(\Gamma)  \normal \theta (\Lambda)$. Let $1\neq g \in G$,
choose $\lambda \in \Lambda$ such that $\nu(\lambda) = g$, and let
$h = \theta(\lambda)$. Then $h = f g$ for some $f\in A^G$.
Then
\[
A\leq A \cap A^{h} = A\cap A^{f g} = A\cap A^g = 1.
\]
This contradiction implies that either $\Gamma = 1$ or
$\Lambda/\Gamma=1$, as desired.
\end{proof}

\begin{proposition}\label{prop:normalprojectivepairs}
Let $(\Gamma,\Lambda)$ be a $\C$-projective pair. If $\Gamma\neq
\Lambda$, then $\bigcap_{x\in \Lambda}\Gamma^x = 1$.
\end{proposition}

\begin{proof}
Let $\Gamma_0 = \bigcap_{x\in \Lambda}\Gamma^x$. By
Corollary~\ref{cor:splittingpp} there exists $\Lambda_0$ such that
$\Gamma_0 = \Lambda_0\cap \Gamma$ and $\Gamma\Lambda_0=\Lambda$. In
particular $(\Lambda_0:\Gamma_0) = (\Lambda:\Gamma)\neq 1$, i.e.\
$\Gamma_0 \neq \Lambda_0$. Moreover, by
Proposition~\ref{prop:ppundersubgroups}, $(\Gamma_0, \Lambda_0)$ is
a $\C$-projective pair. But since $\Gamma_0\normal \Lambda_0$ and
$\Gamma_0\neq \Lambda_0$, Lemma~\ref{lem:normalprojectivepairs}
implies $\Gamma_0=1$.
\end{proof}

\begin{corollary}\label{cor:openpp}
Let $(\Gamma,\Lambda)$ be a $\C$-projective pair. Assume that
$\Gamma$ is open in $\Lambda$. Then $\Gamma=\Lambda$.
\end{corollary}

\begin{proof}
Assume that $\Gamma\neq \Lambda$ (and in particular, $\Lambda\neq
1$). Since $\Gamma$ is open, the normal core $\bigcap_{x\in \Lambda}
\Gamma^x$ is also open. By
Proposition~\ref{prop:normalprojectivepairs}, $\bigcap_{x\in
\Lambda} \Gamma^x=1$. Consequently $\Lambda/\bigcap_{x\in \Lambda}
\Gamma^x = \Lambda$ is a nontrivial finite group. This contradicts
the fact that $\Lambda$ is $\C$-projective, and hence the assertion.
\end{proof}

\begin{proposition}
Let $\Lambda$ be a $\C$-profinite group and $\Gamma$ a $p$-Sylow
subgroup. Assume that $\Lambda$ has a non-abelian simple quotient
that is divisible by $p$. Then the pair $(\Gamma,\Lambda)$ is not
$\C$-projective.
\end{proposition}

\begin{proof}
Assume the contrary, i.e.\ $(\Gamma,\Lambda)$ is $\C$-projective.
Hence, by Corollary~\ref{cor:semidirectproductpp}, $\Lambda =
M\rtimes \Gamma$. Note that $p\nmid (\Lambda:\Gamma) = |M|$ since
$\Gamma$ is $p$-Sylow. Let $\psi\colon \Lambda \to S$ be an
epimorphism onto a non-abelian simple group of order divisible by
$p$. Then $\psi(M)\neq S$. We thus get that $\psi(M)=1$ (since
$\psi(M)\normal \psi (\Lambda)=S$).

On the other hand, $\psi(\Gamma)$ is a proper subgroup of $S$.
(Otherwise $S$ would be a $p$-group, hence solvable.) The assertion
now follows from the contradiction $S = \psi(\Lambda) =
\psi(M)\psi(\Gamma) = \psi(\Gamma)< S$.
\end{proof}

\section{Applications to PAC extensions}
In this section we shall use the following notation from
\cite[Section~2]{Bary-SorokerPACext}. An embedding problem
$(\mu\colon \gal(K) \to A, \alpha \colon G\to A)$ for a field $K$ is
called \textbf{geometric} if there exists a $G$-extension $F/E$ such
that $E$ is regular over $K$ of transcendence degree $1$, if we set
$L = F\cap K_s$, then $L$ is an $A$-extension of $K$, and the
restriction map $\gal(F/K(x))\to \gal(L/K)$ coincides with $\alpha$.
If in addition $E = K(x)$, then the embedding problem is called
\textbf{rational}.

A weak solution of a geometric embedding problem is called
\textbf{geometric} if it is induced from a $K$-rational place $\phi$
of $E$ that is unramified in $F$.

The notion of a double embedding problem for a separable algebraic
field extension $M/K$ comes from the pair $(\gal(M),\gal(K))$. A
double embedding problem is called \textbf{rational} if the higher
embedding problem is rational. A weak solution $(\eta,\theta)$ of a
rational double embedding problem is called \textbf{geometric} if
both $\eta$ and $\theta$ are geometric, say w.r.t.\ $\phi$ and
$\psi$ respectively, and $\psi$ is the restriction of $\phi$ to
$K(x)$ (the field defining the rational higher embedding problem).
See \cite[Sections~3.2 and 3.3]{Bary-SorokerPACext}.

\subsection{Proof of Theorem~\ref{prop:projpair-PACext}}
\begin{trivlist}
\item
Let $M$ be a PAC extension of a PAC field $K$. To prove (a), we need
to show that the pair $(\gal(M),\gal(K))$ is projective. First note
that $\gal(M) \cong \gal(M\cap K_s)$ via the restriction map and
$(M\cap K_s)/K$ is PAC (\cite[Theorem~4.2]{Bary-SorokerPACext}).
Thus we can replace $M$ and $M\cap K_s$, if necessary, to assume
that that $M/K$ is separable and algebraic. Let $\Gamma = \gal(M)$
and $\Lambda=\gal(K)$.

Since $K$ is PAC, $\Lambda$ is projective \cite[Theorem
11.6.2]{FriedJarden2005}. By
Corollary~\ref{cor:ppreductiontodoublesplit}, to show that
$(\Gamma,\Lambda)$ is projective it suffices to solve a split double
embedding problem \eqref{double embedding problem}. Over PAC fields
any finite split embedding problem is rational (see e.g.\
\cite{Pop1996,HaranJarden1998}), and hence any split DEP is
rational. By \cite[Proposition~4.5]{Bary-SorokerPACext} there exists
a weak geometric solution, and in particular a weak solution, of any
finite split DEP.

For (b) assume that $M$ is an algebraic extension of a PAC field $K$
and that $(\gal(M),\gal(K))$ is projective. We have to prove that
$M/K$ is PAC.

Assume that $M/K$ is also separable. We use
\cite[Proposition~4.5]{Bary-SorokerPACext} which says that it
suffices to geometrically solve (in the weak sense) each finite
rational double embedding problem. Since $(\gal(M),\gal(K))$ is
projective, the double embedding problem is weakly solvable. By
\cite[Corollary~3.4]{Bary-SorokerPACext} every weak solution is
geometric, and hence the assertion.

In the general case, let $N = M\cap K_s$. Then $\gal(M) = \gal(N)$
and $N/K$ is separable. Then $N/K$ is PAC. Then $M/K$ is PAC since
$M/N$ is purely inseparable (\cite[Corollary~2.3]{JardenRazon1994}).

For (c) let $(\Gamma,\Lambda)$ be a projective pair. We need to
construct a PAC extension $M/K$ such that $\Gamma = \gal(M)$,
$\Lambda=\gal(K)$.

By \cite[Corollary 23.1.2]{FriedJarden2005}, there exists a PAC
field $K$ such that $\gal(K)\cong \Lambda$ (since $\Lambda$ is
projective). Let $M$ be the fixed field of $\Gamma$, i.e., $\gal(M)
= \Gamma$. Since $(\Gamma,\Lambda)$ is projective, by (b), $M/K$ is
PAC. \qed
\end{trivlist}

Theorem~\ref{prop:projpair-PACext} group theoretically describes the
structure of a PAC extension $M$ of a PAC field $K$. Removing the
condition that $K$ is PAC gives the following more general problem.

\begin{problem}\label{prob:CharPAC}
Describe, purely group theoretically, the pairs $(\gal(M),
\gal(K))$, where $M/K$ is PAC (and $K$ is arbitrary).
\end{problem}

Note that this problem generalizes the classical problem of
characterizing the class of absolute Galois groups out of all the
profinite groups, and hence is much more difficult. Indeed $K_s/K$
is PAC whenever $K$ is infinite. Hence a description of the pair
$(1,\gal(K))$ clearly gives a description of $\gal(K)$.

\subsection{New examples of PAC extensions}
We follow the method outlined in the introduction to construct PAC
extensions.

\begin{proposition}
Let $K_0$ be a field which has a PAC extension $K/K_0$. Assume that
$\gal(K)$ is free of infinite rank $\kappa$. Then for any projective
group $P$ of rank $\leq \kappa$ there exists a PAC extension $M/K_0$
such that $P\cong \gal(M)$.
\end{proposition}

\begin{proof}
By \cite[Theorem~4.2]{Bary-SorokerPACext} we can assume that $K/K_0$
is a separable algebraic extension. By Corollary~\ref{cor:ppHL} $P$
embeds into $\gal(K)$ in such a way that $(P,\gal(K))$ is
projective. Theorem~\ref{prop:projpair-PACext} now implies that for
the fixed field $M$ of $P$ (i.e.\ $\gal(M) = P$), the extension
$M/K$ is PAC. Now the transitivity of PAC extensions
(\cite[Theorem~5]{Bary-SorokerPACext}) implies that $M/K_0$ is PAC.
\end{proof}

Recall that a Galois extension $N/K$ is unbounded if the set
\[\{\ord(\sigma)\mid \sigma\in \gal(N/K)\}\subseteq \bbN\cup \infty\]
is unbounded.

\begin{corollary}\label{cor:ppAbelianExtension}
Let $P$ be a projective group of at most countable rank, let $K_0$
be a countable Hilbertian field, and let $K/K_0$ be an unbounded
abelian extension. Then $K$ has a PAC extension $M$ such that
$P\cong \gal(M)$.
\end{corollary}

\begin{proof}
In the proof of \cite[Proposition 3.8]{Razon1997} it is shown that
there exists a PAC extension $M/K$ such that $\gal(M)\cong
\Fhat_\omega$. Hence the previous proposition implies that there
exists a PAC extension $N/K$ with $\gal(N) \cong P$.
\end{proof}

\begin{remark}
A noteworthy  special case of the last result is when $K_0$ is a
finitely generated infinite field and $K$ is its maximal abelian
extension.
\end{remark}

\begin{corollary}
Let $P$ be a countable projective group. Then there exists a
Hilbertian field $K$ and a PAC extension $M/K$ such that
$\gal(M)\cong P$.
\end{corollary}

\begin{proof}
Let $K_0 = \bbQ$ and $K = \bbQ^{\ab}$. Then $K$ is Hilbertian
\cite[Theorem~16.11.3]{FriedJarden2005} and there exists a PAC
extension $M/K$ such that $\gal(M) = P$
(Corollary~\ref{cor:ppAbelianExtension}).
\end{proof}

% ----------------------------------------------------------------

\end{document}